\documentclass[a4paper, 11pt]{article}
\usepackage{mathrsfs}

\usepackage{epsfig}
\usepackage{amsmath}
\usepackage{amssymb}
\usepackage{latexsym}
\usepackage{amsfonts}
\usepackage{amsthm}
\usepackage[numbers,sort&compress]{natbib}
\usepackage{graphicx}
\ifx\pdfoutput\undefined  \else
 \fi
\usepackage[centerlast]{caption2}
\usepackage{color}

\numberwithin{figure}{section} \numberwithin{definition}{section}
\numberwithin{observation}{section} \numberwithin{lemma}{section}
\numberwithin{theorem}{section} \numberwithin{proposition}{section}
\numberwithin{conjecture}{section}

\begin{document}

\title{
{On the Domination Number of \\ Generalized Petersen Graphs
$P(ck,k)$} \footnote{The research is supported by NSFC (60973014,
11001035), Specialized Research Fund for the Doctoral Program of
Higher Education (200801411073) and Research Foundation of DLUT.}
\author{
Haoli Wang, \ Xirong Xu, \ Yuansheng Yang\footnote {corresponding
author's email : yangys@dlut.edu.cn}, \\
Department of Computer Science \\
Dalian University of Technology, Dalian, 116024, P.R. China\\
\\
Guoqing Wang \\
Center for Combinatorics, LPMC-TJKLC \\
Nankai University, Tianjin, 300071, P.R. China\\
}}

\date{}
\maketitle
\begin{abstract}
\noindent Let $G=(V(G),E(G))$ be a simple connected and undirected
graph with vertex set $V(G)$ and edge set $E(G)$. A set $S \subseteq
V(G)$ is a $dominating$ $set$ if for each $v \in V(G)$ either $v \in
S$ or $v$ is adjacent to some $w \in S$. That is, $S$ is a
dominating set if and only if $N[S]=V(G)$. The domination number
$\gamma(G)$ is the minimum cardinalities of minimal dominating sets.
In this paper, we give an improved upper bound on the domination
number of generalized Petersen graphs $P(ck,k)$ for $c\geq 3$ and
$k\geq 3$. We also prove that $\gamma(P(4k,k))=2k+1$ for even $k$,
$\gamma(P(5k,k))=3k$ for all $k\geq 1$, and
$\gamma(P(6k,k))=\lceil\frac{10k}{3}\rceil$ for $k\geq 1$ and $k\neq
2$.



\end{abstract}

\section{Introduction}

\ \ \ \ Let $G=(V(G),E(G))$ be a simple connected and undirected
graph with vertex set $V(G)$ and edge set $E(G)$. The open
neighborhood and the closed neighborhood of a vertex $v\in V(G)$ are
denoted by $N(v)=\{u \in V(G)\ :\ vu \in E(G)\}$ and $N[v]=N(v)\cup
\{v\}$, respectively. For a vertex set $S \subseteq V(G)$,
$N(S)=\underset{v \in S}\cup N(v)$ and $N[S]=\underset{v \in S}\cup
N[v]$. For $S\subseteq V(G)$, let $\langle S\rangle$ be the subgraph
induced by $S$.

A set $S \subseteq V(G)$ is a $dominating$ $set$ if for each $v \in
V(G)$ either $v \in S$ or $v$ is adjacent to some $w \in S$. That
is, $S$ is a dominating set if and only if $N[S]=V(G)$. The
domination number of $G$, denoted by $\gamma(G)$, is the minimum
cardinalities of minimal dominating sets. A subset $S\subset V(G)$
is efficient dominating set or a perfect dominating set if each
vertex of $G$ is dominated by exactly one vertex in $D$. For a more
detailed treatment of domination-related parameters and for
terminology not defined here, the reader is referred to
\cite{HHS98}.

In recent years, domination and its variations on the class of
generalized Petersen graph have been studied extensively
\cite{BBP08,EJM09,FGS,XYB07,XYB08,YKX,Z02,ZZW}. The generalized
Petersen graph $P(n,k)$ is defined to be a graph on $2n$ vertices
with $V(P(n,k)) = \{v_i,u_{i}:0 \leq i \leq n-1\}$ and $E(P(n,k)) =
\{v_iv_{i+1},v_{i}u_{i},u_{i}u_{i+k}:0 \leq i \leq n-1$, subscripts
are taken modulo $n \}$. In 2009, B. Javad Ebrahimi et al
\cite{EJM09} proved a necessary and sufficient condition for the
generalized Petersen graphs to have an efficient dominating set.

\noindent \textbf{Lemma 1.1. \cite{EJM09}} If $P(n,k)$ has an
efficient dominating set, then $\gamma(P(n,k))$ $=\frac{n}{2}$ and
$n\equiv 0$ (mod 4).

\noindent \textbf{Theorem 1.2. \cite{EJM09}} A generalized Petersen
graph $P(n,k)$ has an efficient dominating set if and only if
$n\equiv 0$ (mod 4) and $k$ is odd.


Recently, Weiliang Zhao et al \cite{ZZW} have started to study the
domination number of the generalized Petersen graphs $P(ck,k)$,
where $c\geq 3$ is a constant. They obtained upper bound on
$\gamma(P(ck,k))$ for $c\geq 3$ as follows:
$$\begin{array}{llll}
\gamma(P(ck,k))\leq& \left \{\begin{array}{llll}
               \frac{c}{3}\lceil\frac{5k}{3}\rceil, & \mbox{ if } c\equiv 0 \mbox{ (mod } 3);\\
               \lceil\frac{c}{3}\rceil\lceil\frac{5k}{3}\rceil-\lceil\frac{2k}{3}\rceil, & \mbox{ if } c\equiv 1 \mbox{ (mod } 3);\\
               \lceil\frac{c}{3}\rceil\lceil\frac{5k}{3}\rceil-\lceil\frac{2k}{3}\rceil+\lceil\frac{k}{3}\rceil, & \mbox{ if } c\equiv 2 \mbox{ (mod } 3).\\

              \end{array}
           \right . \\
\end{array}$$
They also determined the domination number of $P(3k,k)$ for $k\geq
1$ and the domination number of $P(4k,k)$ for odd $k$.

In this paper, we study the domination number of generalized
Petersen graphs $P(ck,k)$. We give an improved upper bound on the
domination number of $P(ck,k)$ for $c\geq 3$ and $k\geq 3$. We also
prove that $\gamma(P(4k,k))=2k+1$ for even $k$, $\gamma(P(5k,k))=3k$
for all $k\geq 1$, and $\gamma(P(6k,k))=\lceil\frac{10k}{3}\rceil$
for $k\geq 1$ and $k\neq 2$.

Throughout the paper, the subscripts are taken modulo $n$ when it is
unambiguous.

\section{General upper bound of $P(ck,k)$}

\ \ \ \ In this section, we shall give an improved upper bound on
the domination number of $P(ck,k)$ for general $c$.

\noindent \textbf{Theorem 2.1.} For any constant $c\geq 3$ and
$k\geq 3$,
 {\footnotesize
$$\begin{array}{llll}
\gamma(P(ck,k))\leq& \left \{\begin{array}{llll}
               \frac{ck}{2}+\alpha, & \mbox{ if } c\equiv 0 \mbox{ (mod } 4);\\
               \frac{ck}{2}+\frac{k}{2}-1+\alpha, & \mbox{ if } c\equiv 1,2 \mbox{ (mod } 4) \mbox{ and } k\equiv 0 \mbox{ (mod } 2);\\
               \frac{ck-1}{2}+\frac{k+1}{2}+\alpha, & \mbox{ if } c\equiv 1 \mbox{ (mod } 4) \mbox{ and } k\equiv 1 \mbox{ (mod } 2);\\
               \frac{ck}{2}+\frac{k+1}{2}+\alpha, & \mbox{ if } c\equiv 2 \mbox{ (mod } 4) \mbox{ and } k\equiv 1 \mbox{ (mod } 4);\\
               \frac{ck}{2}+\frac{k-1}{2}+\alpha, & \mbox{ if } c\equiv 2 \mbox{ (mod } 4) \mbox{ and } k\equiv 3 \mbox{ (mod } 4);\\
               \lfloor\frac{ck}{2}\rfloor+\lfloor\frac{k}{4}\rfloor+1+\alpha, & \mbox{ if } c\equiv 3 \mbox{ (mod } 4) \mbox{ and } k\neq 4,8;\\
               \frac{ck}{2}+\frac{k}{4}+\alpha, & \mbox{ if } c\equiv 3 \mbox{ (mod } 4) \mbox{ and } k=4,8;\\
              \end{array}
           \right . \\
\end{array}$$}
where {\footnotesize
$$\begin{array}{llll}
\alpha=& \left \{\begin{array}{llll}
               0, & \mbox{ if } k\equiv 1 \mbox{ (mod } 2);\\
               \lfloor\frac{c}{4}\rfloor, & \mbox{ if } k\equiv 0 \mbox{ (mod } 2).\\
               \end{array}
           \right . \\
\end{array}$$
}
\begin{proof} To show this upper bound, it suffices to give a dominating set $S$ with the cardinality equaling
to the values mentioned in this theorem. Let $n=ck$,
$m=\lfloor\frac{n}{4}\rfloor$ and $t=n$ mod 4. Then $n=4m+t$.

For $k\equiv1 \pmod 2$, let $S_0=A\cup B$, where {\footnotesize
$$A=\{v_{4i}:0\leq i\leq m-1\}\mbox{ \ \ and \ \ }
B=\{u_{4i+2}:0\leq i\leq m-1\},$$} and let {\footnotesize
$$\begin{array}{llll}
S=\left \{\begin{array}{llll}
               S_0,  \ \ \ \ \ \ \ \ \ \ \ \ \ \ \ \ \ \ \ \ \ \ \ \ \ \ \ \ \ \ \ \ \ \ \ \ \ \ \ \ \ \ \ \ \ \ \ \ \ \mbox{ if } c\equiv 0 \mbox{ (mod } 4);\\
               S_0\cup\{u_{n-2-4i},u_{n-4-4i}:0\leq i\leq \lfloor \frac{k}{4}\rfloor-1\}\cup\{u_{n-1}\}, \\
               \ \ \ \ \ \ \ \ \ \ \ \ \ \ \ \ \ \ \ \ \ \ \ \ \ \ \ \ \ \ \ \ \ \ \ \ \ \ \ \ \ \ \ \ \ \ \ \ \ \ \ \ \ \ \mbox{ if } c\equiv 1 \mbox{ (mod } 4) \mbox{ and } k\equiv 1 \mbox{ (mod } 4);\\
               S_0\cup\{u_{n-2-4i},u_{n-4-4i}:0\leq i\leq \lceil \frac{k}{4}\rceil-1\}\cup\{v_{n-3}\}, \\
               \ \ \ \ \ \ \ \ \ \ \ \ \ \ \ \ \ \ \ \ \ \ \ \ \ \ \ \ \ \ \ \ \ \ \ \ \ \ \ \ \ \ \ \ \ \ \ \ \ \ \ \ \ \ \mbox{ if } c\equiv 1 \mbox{ (mod } 4) \mbox{ and } k\equiv 3 \mbox{ (mod } 4);\\
               S_0\cup\{u_{n-2-4i},u_{n-5-4i}:0\leq i\leq \lfloor \frac{k}{4}\rfloor-1\}\cup\{u_{n-1},u_{n-3}\}, \\
               \ \ \ \ \ \ \ \ \ \ \ \ \ \ \ \ \ \ \ \ \ \ \ \ \ \ \ \ \ \ \ \ \ \ \ \ \ \ \ \ \ \ \ \ \ \ \ \ \ \ \ \ \ \ \mbox{ if } c\equiv 2 \mbox{ (mod } 4) \mbox{ and } k\equiv 1 \mbox{ (mod } 4);\\
               S_0\cup\{u_{n-2-4i},u_{n-3-4i}:0\leq i\leq \lceil\frac{k}{4}\rceil-1\},  \\
               \ \ \ \ \ \ \ \ \ \ \ \ \ \ \ \ \ \ \ \ \ \ \ \ \ \ \ \ \ \ \ \ \ \ \ \ \ \ \ \ \ \ \ \ \ \ \ \ \ \ \ \ \ \ \mbox{ if } c\equiv 2 \mbox{ (mod } 4) \mbox{ and } k\equiv 3 \mbox{ (mod } 4);\\
               S_0\cup\{u_{n-2-4i}:0\leq i\leq \lfloor \frac{k}{4}\rfloor\}\cup\{v_{n-3}\},  \ \mbox{ if } c\equiv 3 \mbox{ (mod } 4) \mbox{ and } k\equiv 1 \mbox{ (mod } 4);\\
               S_0\cup\{u_{n-2-4i}:0\leq i\leq \lfloor \frac{k}{4}\rfloor\},  \ \ \ \ \ \ \ \ \ \ \ \ \ \ \mbox{ if } c\equiv 3 \mbox{ (mod } 4) \mbox{ and } k\equiv 3 \mbox{ (mod } 4).\\
              \end{array}
           \right . \\
\end{array}$$}
It is not hard to check that {\footnotesize
$$\begin{array}{llll}
|S|=\left \{\begin{array}{llll}
               \frac{ck}{2}, \mbox{if } c\equiv 0 \mbox{ (mod } 4);\\
               2\times\lfloor\frac{ck}{4}\rfloor+2\times\lfloor\frac{k}{4}\rfloor+1=\frac{ck-1}{2}+\frac{k+1}{2}, \ \ \mbox{if } c\equiv 1 \mbox{ (mod } 4)\mbox{ and } k\equiv 1 \mbox{ (mod } 4);\\
               2\times\lfloor\frac{ck}{4}\rfloor+2\times\lceil\frac{k}{4}\rceil+1=\frac{ck-1}{2}+\frac{k+1}{2}, \ \ \mbox{if } c\equiv 1 \mbox{ (mod } 4)\mbox{ and } k\equiv 3 \mbox{ (mod } 4);\\
               2\times\lfloor\frac{ck}{4}\rfloor+2\times\lfloor\frac{k}{4}\rfloor+2=\frac{ck}{2}+\frac{k+1}{2}, \ \ \ \ \ \mbox{if } c\equiv 2 \mbox{ (mod } 4) \mbox{ and } k\equiv 1 \mbox{ (mod } 4);\\
               2\times\lfloor\frac{ck}{4}\rfloor+2\times\lceil\frac{k}{4}\rceil=\frac{ck}{2}+\frac{k-1}{2}, \ \ \ \ \ \ \ \ \ \ \mbox{if } c\equiv 2 \mbox{ (mod } 4) \mbox{ and } k\equiv 3 \mbox{ (mod } 4);\\
               2\times\lfloor\frac{ck}{4}\rfloor+\lfloor\frac{k}{4}\rfloor+2=\lfloor\frac{ck}{2}\rfloor+\lfloor\frac{k}{4}\rfloor+1, \ \ \ \mbox{if } c\equiv 3 \mbox{ (mod } 4) \mbox{ and } k\equiv 1 \mbox{ (mod } 4);\\
               2\times\lfloor\frac{ck}{4}\rfloor+\lfloor\frac{k}{4}\rfloor+1=\lfloor\frac{ck}{2}\rfloor+\lfloor\frac{k}{4}\rfloor+1, \ \ \ \mbox{if } c\equiv 3 \mbox{ (mod } 4) \mbox{ and } k\equiv 3 \mbox{ (mod } 4).\\
              \end{array}
           \right . \\
\end{array}$$}

For $k\equiv0 \pmod 2$, let $m_2=\lfloor\frac{c}{4}\rfloor$ and
$r=c$ mod 4. Denote

\noindent $S_{40}=A_{40}\cup B_{40}\cup C_{40}\cup D_{40}\cup
E_{40}$, where
{\footnotesize $$\begin{array}{llll}
        A_{40}=\{v_{4kj+2+4i},u_{4kj+4i}:0\leq i\leq \frac{k}{4}-1, \ \ 0\leq j\leq m_2-1\}, \\
        B_{40}=\{v_{4kj+k+1+4i},u_{4kj+k+3+4i}:0\leq i\leq \frac{k}{4}-1, \ \ 0\leq j\leq m_2-1\}, \\
        C_{40}=\{v_{4kj+2k+4i},u_{4kj+2k+2+4i}:0\leq i\leq \frac{k}{4}-1, \ \ 0\leq j\leq m_2-1\}, \\
        D_{40}=\{v_{4kj+3k+3+4i},u_{4kj+3k+1+4i}:0\leq i\leq \frac{k}{4}-1, \ \ 0\leq j\leq m_2-1\}, \\
        E_{40}=\{v_{4kj+3k}:0\leq j\leq m_2-1\}, \\
\end{array}$$
} \noindent $S_{42}=A_{42}\cup B_{42}\cup C_{42}\cup D_{42}\cup
E_{42}$, where
{\footnotesize $$\begin{array}{llll}
        \ \ \ A_{42}=\{v_{4kj+4i},u_{4kj+2+4i}:0\leq i\leq \frac{k-2}{4}-1, \ \ 0\leq j\leq m_2-1\}, \\
        \ \ \ B_{42}=\{v_{4kj+k+1+4i},u_{4kj+k-1+4i}:0\leq i\leq \frac{k-2}{4}-1, \ \ 0\leq j\leq m_2-1\}, \\
        \ \ \ C_{42}=\{v_{4kj+2k+2+4i},u_{4kj+2k+4i}:0\leq i\leq \frac{k-2}{4}-1, \ \ 0\leq j\leq m_2-1\}, \\
        \ \ \ D_{42}=\{v_{4kj+3k-1+4i},u_{4kj+3k+1+4i}:0\leq i\leq \frac{k-2}{4}-1, \ \ 0\leq j\leq m_2-1\}, \\
        \ \ \ E_{42}=\{v_{4kj+k-2},v_{4kj+2k-2},v_{4kj+4k-3},u_{4kj+2k-3},u_{4kj+4k-2},:0\leq j\leq m_2-1\}, \\
\end{array}$$
}
and
{\footnotesize$$\begin{array}{llll} S_4=& \left
\{\begin{array}{llll}
               S_{40}, & \mbox{ if } k\equiv 0 \mbox{ (mod }4);\\
               S_{42}, & \mbox{ if } k\equiv 2 \mbox{ (mod }4).\\
              \end{array}
           \right . \\
\end{array}$$}
Then {\footnotesize
$$\begin{array}{llll}
|S_4|=\left \{\begin{array}{llll}
               2\times \frac{k}{4}\times \frac{c-r}{4}\times 4+\frac{c-r}{4}=\frac{(c-r)k}{2}+\frac{c-r}{4}, & \mbox{ if } k\equiv 0 \mbox{ (mod }4);\\
               2\times \frac{k-2}{4}\times \frac{c-r}{4}\times 4+5\times \frac{c-r}{4}=\frac{(c-r)k}{2}+\frac{c-r}{4}, & \mbox{ if } k\equiv 2 \mbox{ (mod }4).\\
              \end{array}
           \right . \\
\end{array}$$
} \ \ \ \ If $c\equiv 0 \pmod 4$, let $S=S_4$. Then
$|S|=\frac{ck}{2}+\frac{c}{4}$.

If $c\equiv 1 \pmod 4$, let {\footnotesize
$$\begin{array}{llll}
S=\left \{\begin{array}{llll}
               S_4\cup \{u_i:n-k+1\leq i\leq n-1\}, & \mbox{if } k\equiv 0 \mbox{ (mod }4);\\
               S_4\cup \{u_i:n-k+1\leq i\leq n-4\} \cup \{v_{n-k},v_{n-3},u_{n-1}\}, & \mbox{if } k\equiv 2 \mbox{ (mod }4).\\
              \end{array}
           \right . \\
\end{array}$$}
Then {\footnotesize
$$\begin{array}{llll}
|S|=\left \{\begin{array}{llll}
               \frac{(c-1)\times k}{2}+\frac{c-1}{4}+k-1=\frac{ck}{2}+\frac{k}{2}+\lfloor\frac{c}{4}\rfloor-1, & \mbox{ if } k\equiv 0 \mbox{ (mod }4);\\
               \frac{(c-1)\times k}{2}+\frac{c-1}{4}+k-4+3=\frac{ck}{2}+\frac{k}{2}+\lfloor\frac{c}{4}\rfloor-1, & \mbox{ if } k\equiv 2 \mbox{ (mod }4).\\
              \end{array}
           \right . \\
\end{array}$$}
\indent If  $c\equiv 2 \pmod 4$, let {\footnotesize
$$\begin{array}{llll}
S=\left \{\begin{array}{llll}
               S_4\cup \{v_{n-2k+2+4i},u_{n-2k+4i}:0\leq i\leq \frac{k}{4}-1\}\\
               \ \ \ \ \cup \{u_i:n-k\leq i\leq n-1\}\setminus \{u_{n-2k}\}, & \mbox{ if } k\equiv 0 \mbox{ (mod }4);\\
               S_4\cup \{v_{n-2k+4i},u_{n-2k+2+4i}:0\leq i\leq \frac{k-2}{4}-1\}\\
               \ \ \ \ \cup \{u_i:n-k-3\leq i\leq n-5\}\cup \{v_{n-3}\}, & \mbox{ if } k\equiv 2 \mbox{ (mod }4).\\
              \end{array}
           \right . \\
\end{array}$$}
Then {\footnotesize
$$\begin{array}{llll}
|S|=\left \{\begin{array}{llll}
               \frac{(c-2)\times k}{2}+\frac{c-2}{4}+2\times \frac{k}{4}+k-1=\frac{ck}{2}+\frac{k}{2}+\lfloor\frac{c}{4}\rfloor-1, \ \ \ \ \ \ \ \ \ \mbox{if } k\equiv 0 \mbox{ (mod }4);\\
               \frac{(c-2)\times k}{2}+\frac{c-2}{4}+2\times \frac{k-2}{4}+k-1+1=\frac{ck}{2}+\frac{k}{2}+\lfloor\frac{c}{4}\rfloor-1, \ \mbox{if } k\equiv 2 \mbox{ (mod }4).\\
              \end{array}
           \right . \\
\end{array}$$}
\indent If $c\equiv 3 \pmod 4$, let {\footnotesize
$$\begin{array}{llll}
S=\left \{\begin{array}{llll}
               S_4 \cup \{v_{n-ik},v_{n-ik+3}:1\leq i\leq 3\} \\
               \ \ \ \ \cup \{u_{n-2k+2},u_{n-k+1}\}\setminus\{v_{n-k}\}, & \mbox{if } k=4;\\
               \\
               S_4 \cup \{v_{n-ik},v_{n-ik+3},v_{n-ik+6}:1\leq i\leq 3\} \\
               \ \ \ \ \cup \{u_{n-3k+4},u_{n-2k+2},u_{n-2k+7},u_{n-k+1},u_{n-k+5}\}, & \mbox{if } k=8;\\
               \\
               S_4\cup \{v_{n-3k+6+4i},u_{n-3k+4+4i}:0\leq i\leq \frac{k}{4}-2\} \\
               \ \ \ \ \cup \{v_{n-2k+9+4i},u_{n-2k+11+4i}:0\leq i\leq \frac{k}{4}-3\} \\
               \ \ \ \ \cup \{v_{n-k+8+4i},u_{n-k+9+4i},u_{n-k+10+4i}:0\leq i\leq \frac{k}{4}-3\} \\
               \ \ \ \ \cup \{v_{n-ik},v_{n-ik+3}:1\leq i\leq 3\} \\
               \ \ \ \ \cup \{v_{n-2k+6},v_{n-k+6},v_{n-1}\} \\
               \ \ \ \ \cup \{u_{n-2k+2},u_{n-2k+7},u_{n-k+1},u_{n-k+5}\}, & \mbox{if } k\equiv 0 \mbox{ (mod }4) \\
               & \mbox{and } k\neq 4,8;\\
               S_4\cup \{v_{n-3k+4i},u_{n-3k+2+4i}:0\leq i\leq \frac{k-2}{4}-1\} \\
               \ \ \ \ \cup \{v_{n-2k+1+4i},u_{n-2k-1+4i}:0\leq i\leq \frac{k-2}{4}\} \\
               \ \ \ \ \cup \{v_{n-k+3+4i},u_{n-k+4i},u_{n-k+1+4i}:0\leq i\leq \frac{k-2}{4}-1\} \\
               \ \ \ \ \cup \{v_{n-2k-2},u_{n-2}\}, & \mbox{if } k\equiv 2 \mbox{ (mod }4). \\
              \end{array}
           \right . \\
\end{array}$$}
Then {\footnotesize
$$\begin{array}{llll}
|S|=\left \{\begin{array}{llll}
               \frac{(c-3)\times k}{2}+\frac{c-3}{4}+8-1=\frac{ck}{2}+\frac{k}{4}+\lfloor\frac{c}{4}\rfloor,  \ \ \ \ \ \ \mbox{if } k=4;\\
               \frac{(c-3)\times k}{2}+\frac{c-3}{4}+14=\frac{ck}{2}+\frac{k}{4}+\lfloor\frac{c}{4}\rfloor,  \ \ \ \ \ \ \ \ \ \ \mbox{if } k=8;\\
               \frac{(c-3)\times k}{2}+\frac{c-3}{4}+2\times (\frac{k}{4}-1)+5\times (\frac{k}{4}-2)+13=\frac{ck}{2}+\frac{k}{4}+\lceil\frac{c}{4}\rceil,  \\
               \ \ \ \ \ \ \ \ \ \ \ \ \ \ \ \ \ \ \ \ \ \ \ \ \ \ \ \ \ \ \ \ \ \ \ \ \ \ \ \ \ \ \ \ \ \ \ \ \ \ \ \ \ \ \ \ \ \ \ \ \mbox{if } k\equiv 0 \mbox{ (mod }4) \mbox{ and } k\neq 4,8;\\
               \frac{(c-3)\times k}{2}+\frac{c-3}{4}+2\times (\frac{k-2}{4}+1)+5\times \frac{k-2}{4}+2=\frac{ck}{2}+\frac{k-2}{4}+\lceil\frac{c}{4}\rceil,  \\
               \ \ \ \ \ \ \ \ \ \ \ \ \ \ \ \ \ \ \ \ \ \ \ \ \ \ \ \ \ \ \ \ \ \ \ \ \ \ \ \ \ \ \ \ \ \ \ \ \ \ \ \ \ \ \ \ \ \ \ \ \mbox{if } k\equiv 2 \mbox{ (mod }4). \\
              \end{array}
           \right. \\
\end{array}$$}
\indent It is not hard to verify that $S$ is a dominating set of
$P(ck,k)$ with cardinality equaling to the values mentioned in this
theorem.
\end{proof}

In Figure 2.1 and Figure 2.2, we show the dominating sets of
$P(ck,k)$ for $3\leq k\leq 10$ and $4\leq c\leq 7$, where the
vertices of dominating sets are in dark.
\begin{figure}[ht]
\centering
\includegraphics[scale=0.5]{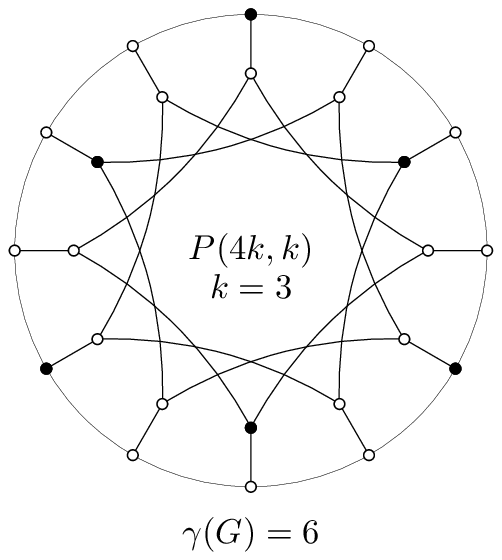} \hspace{2pt}
\includegraphics[scale=0.5]{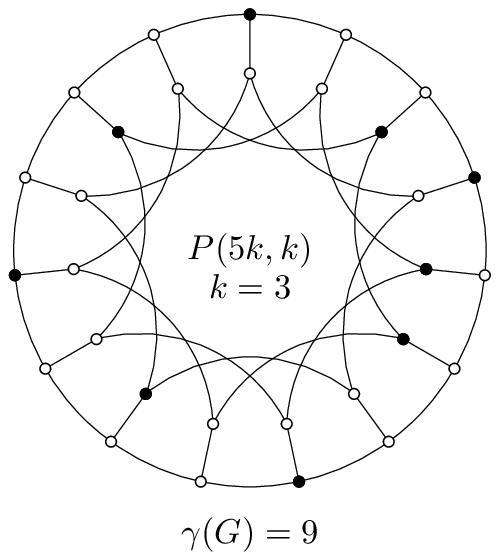} \hspace{2pt}
\includegraphics[scale=0.5]{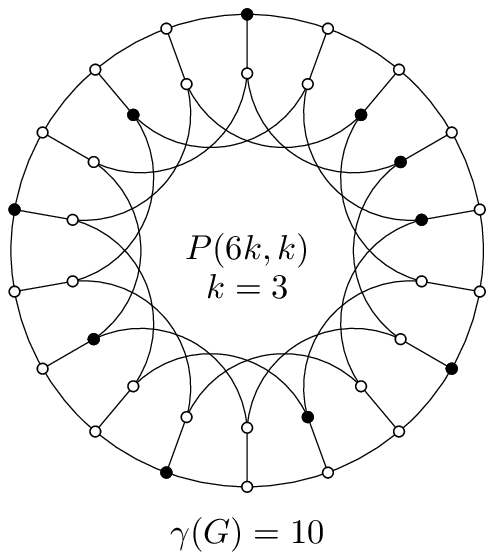} \hspace{2pt}
\includegraphics[scale=0.5]{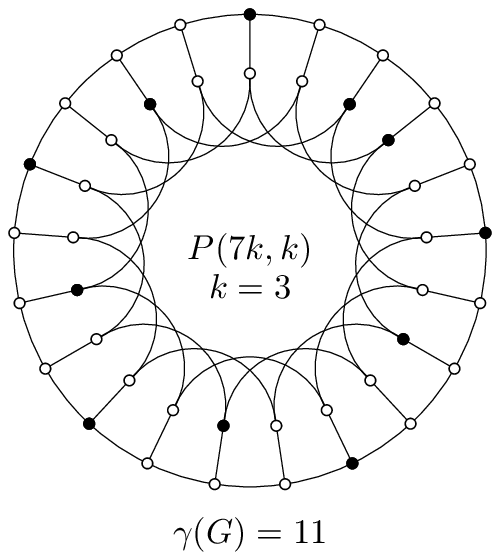} \hspace{2pt}
\includegraphics[scale=0.5]{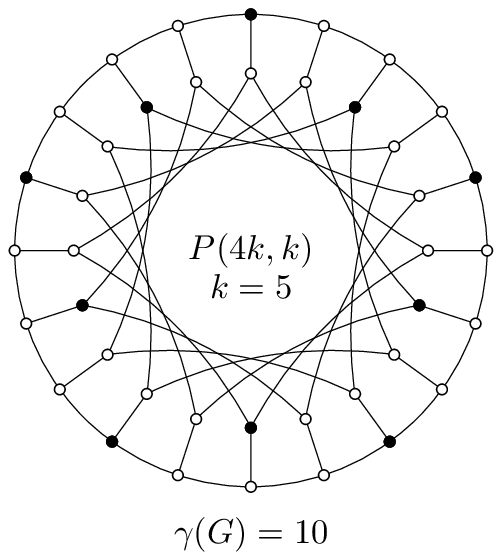} \hspace{2pt}
\includegraphics[scale=0.5]{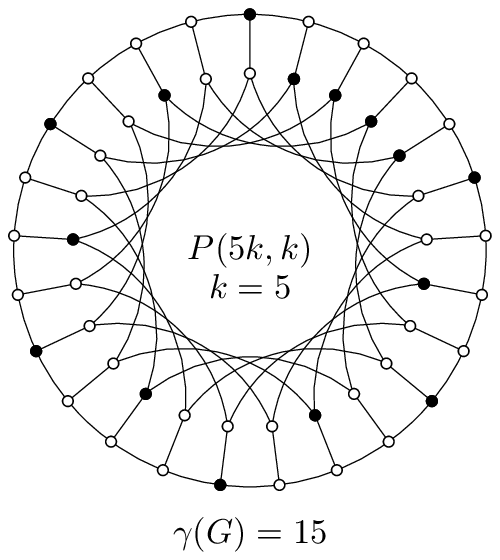} \hspace{2pt}
\includegraphics[scale=0.5]{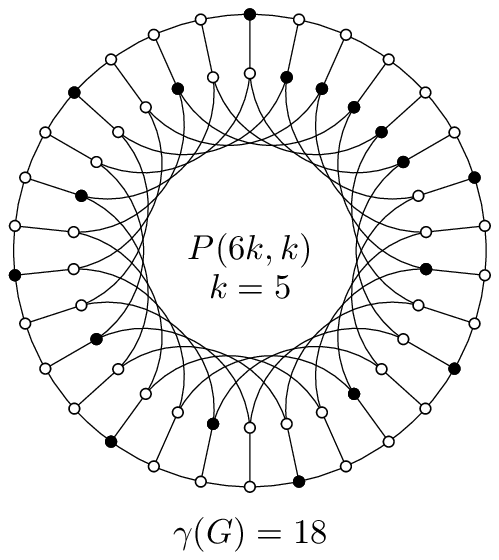} \hspace{2pt}
\includegraphics[scale=0.5]{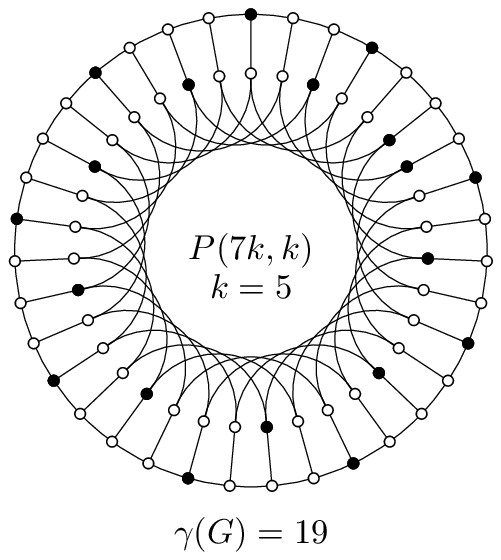} \hspace{2pt}
\includegraphics[scale=0.5]{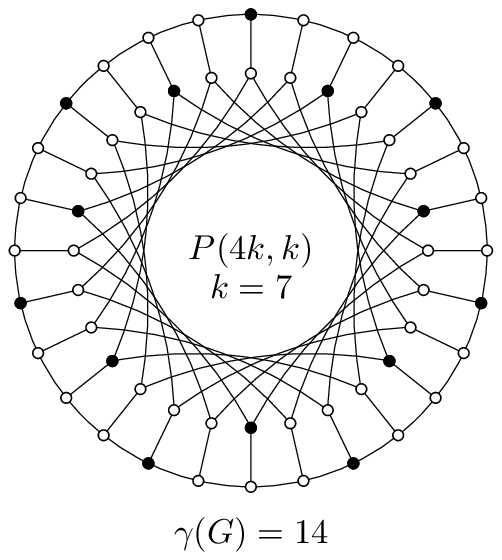} \hspace{2pt}
\includegraphics[scale=0.5]{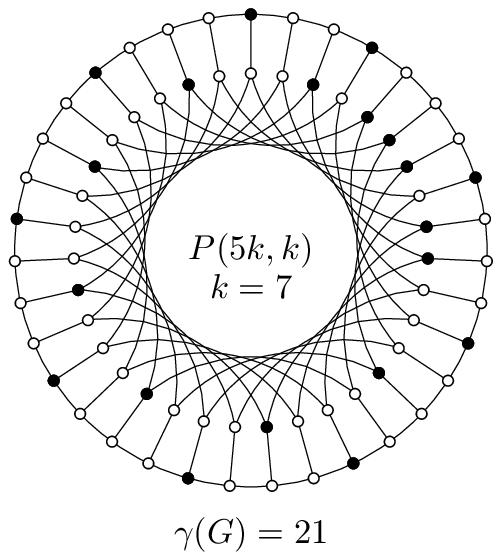} \hspace{2pt}
\includegraphics[scale=0.5]{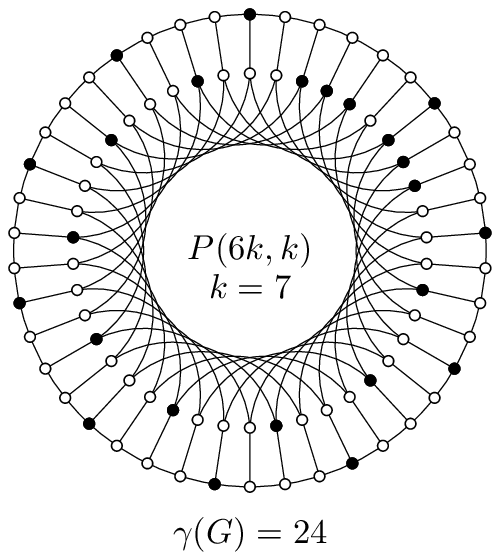} \hspace{2pt}
\includegraphics[scale=0.5]{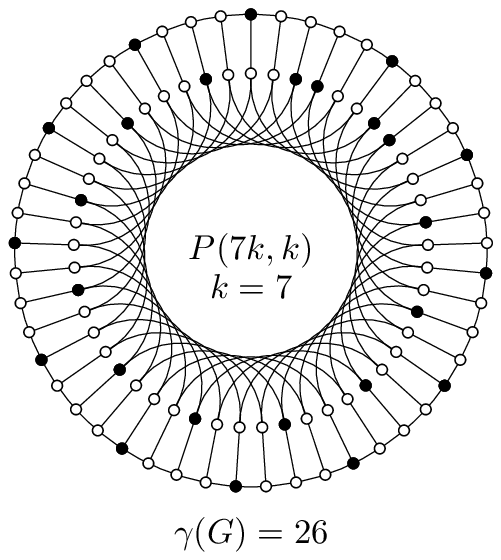} \hspace{2pt}
\includegraphics[scale=0.5]{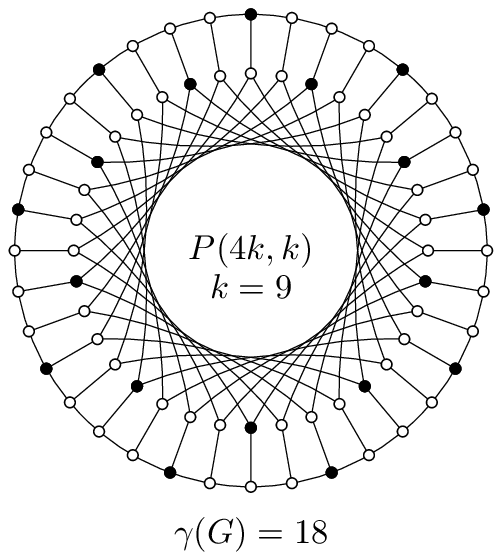} \hspace{2pt}
\includegraphics[scale=0.5]{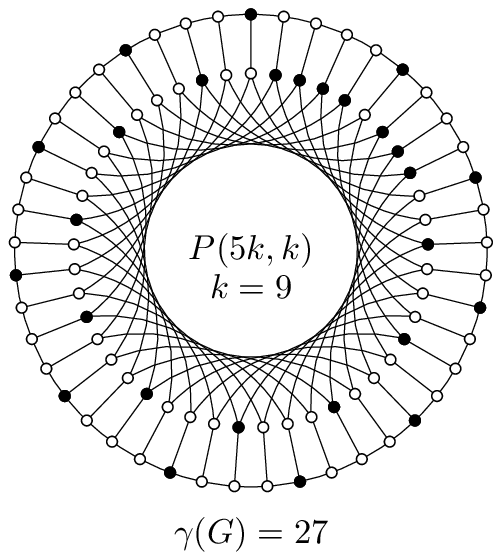} \hspace{2pt}
\includegraphics[scale=0.5]{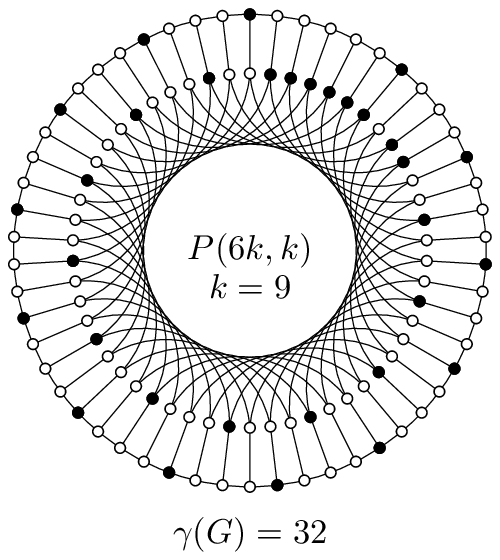} \hspace{2pt}
\includegraphics[scale=0.5]{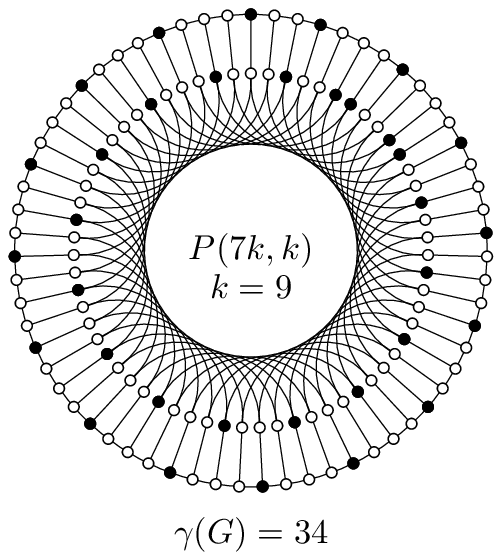} \hspace{2pt}

\centering \vspace{2pt}\small{Figure 2.1:  The dominating sets of
$P(ck,k)$ for $k=3,5,7,9$ and $c=4,5,6,7$}
\end{figure}

As an immediate consequence of Lemma 1.1, Theorem 1.2 and Theorem
2.1, we have the following

{\noindent \bf Theorem 2.2.} For $k\geq 1$,
$$\begin{array}{llll}
\gamma(P(4k,k))=& \left \{\begin{array}{llll}
               2k, & \mbox{ if } k\equiv 1 \mbox{ (mod 2)};\\
               2k+1, & \mbox{ if } k\equiv 0 \mbox{ (mod 2)}.\\
              \end{array}
           \right . \\
\end{array}$$
\begin{figure}[ht]
\centering
\includegraphics[scale=0.5]{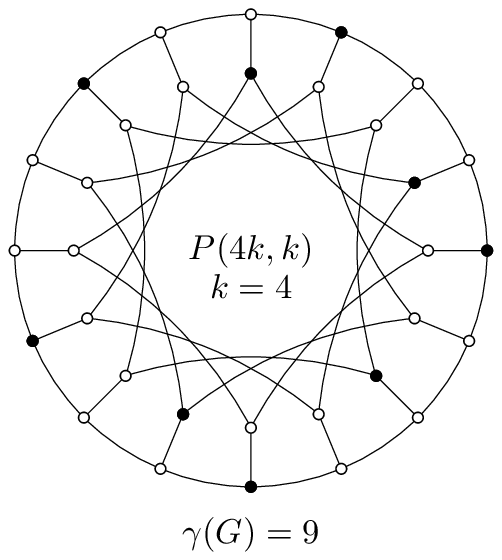} \hspace{2pt}
\includegraphics[scale=0.5]{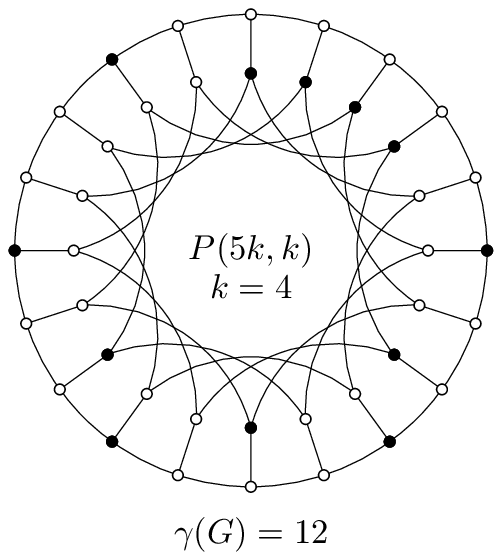} \hspace{2pt}
\includegraphics[scale=0.5]{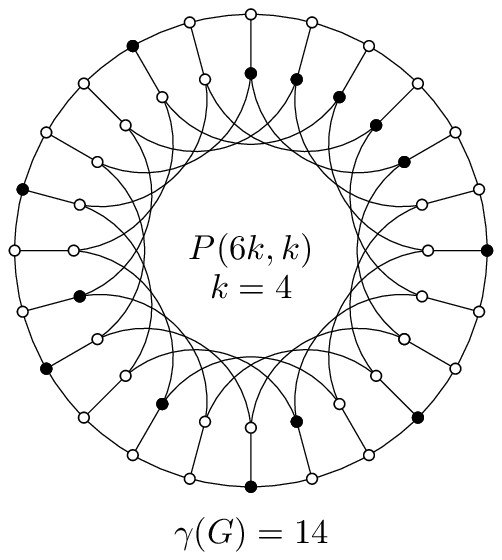} \hspace{2pt}
\includegraphics[scale=0.5]{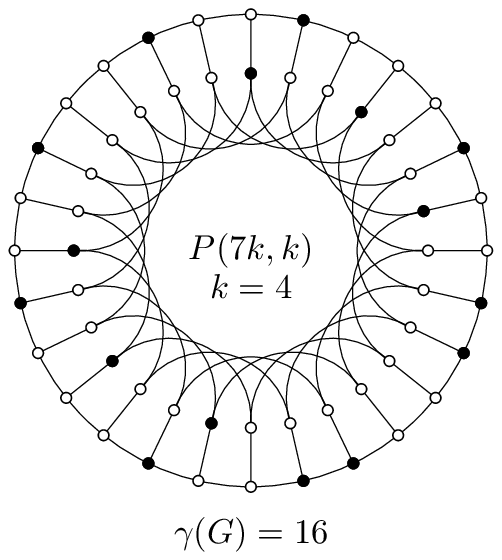} \hspace{2pt}
\includegraphics[scale=0.5]{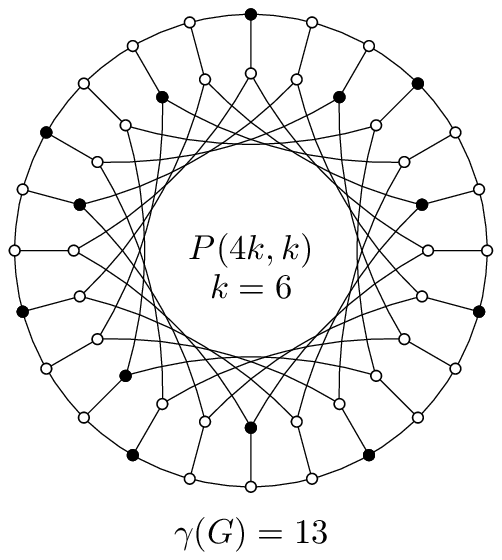} \hspace{2pt}
\includegraphics[scale=0.5]{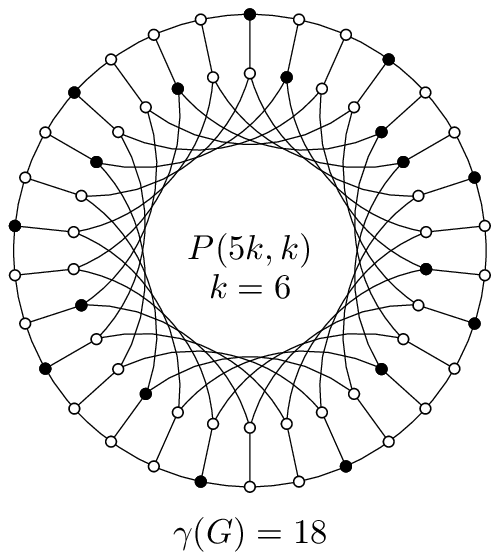} \hspace{2pt}
\includegraphics[scale=0.5]{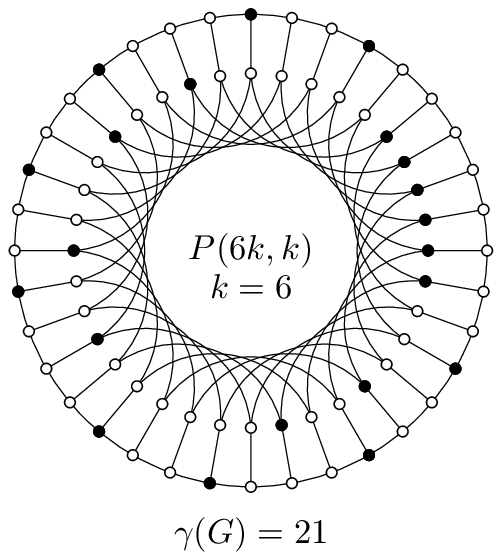} \hspace{2pt}
\includegraphics[scale=0.5]{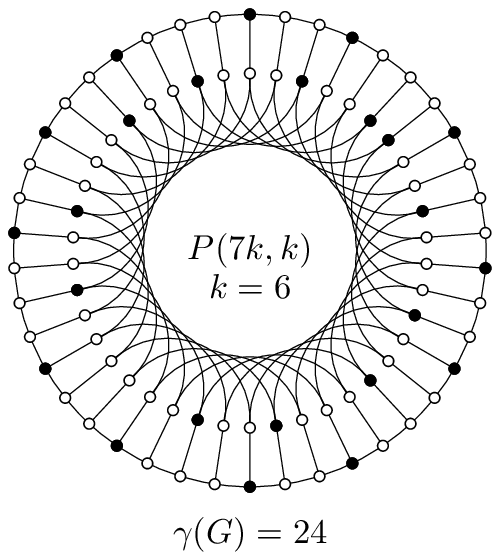} \hspace{2pt}
\includegraphics[scale=0.5]{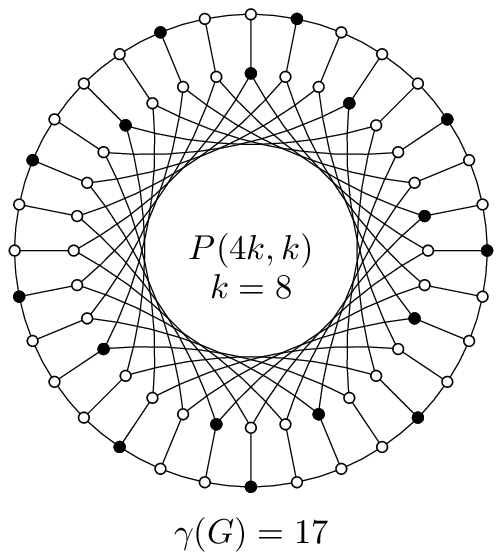} \hspace{2pt}
\includegraphics[scale=0.5]{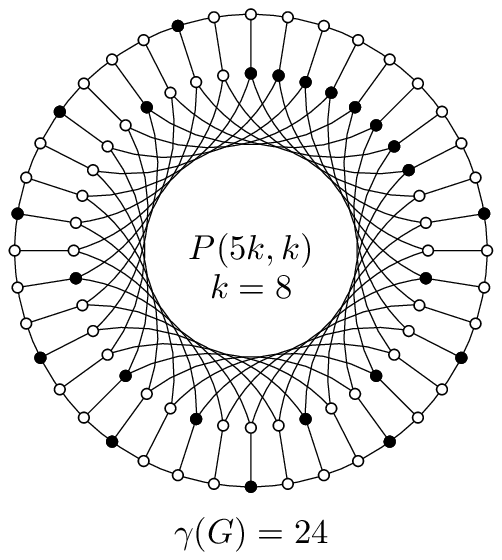} \hspace{2pt}
\includegraphics[scale=0.5]{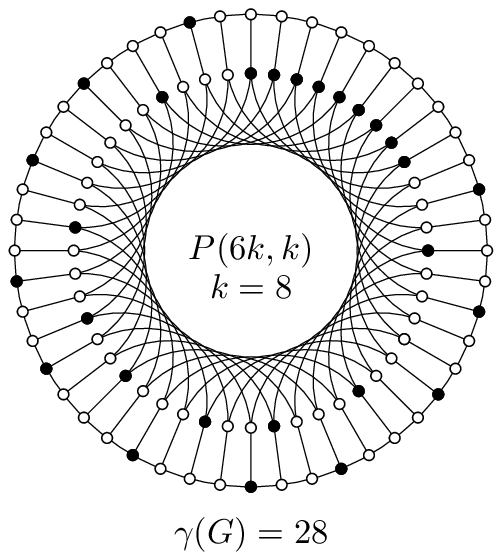} \hspace{2pt}
\includegraphics[scale=0.5]{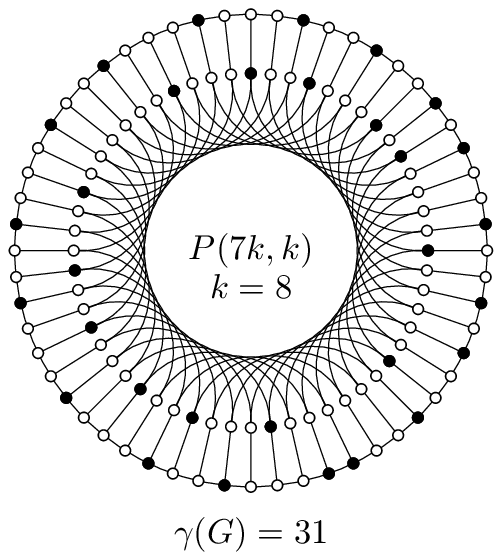} \hspace{2pt}
\includegraphics[scale=0.5]{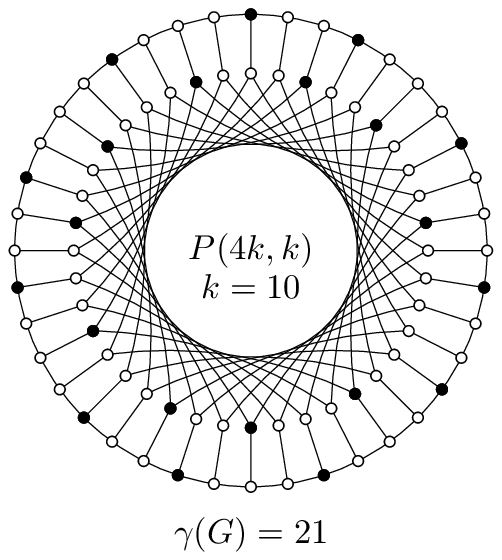} \hspace{2pt}
\includegraphics[scale=0.5]{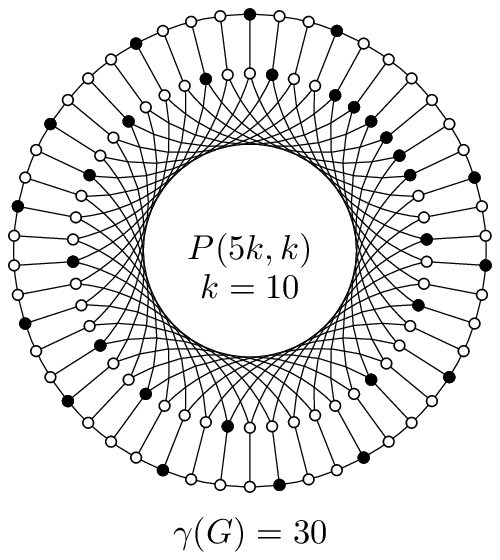} \hspace{2pt}
\includegraphics[scale=0.5]{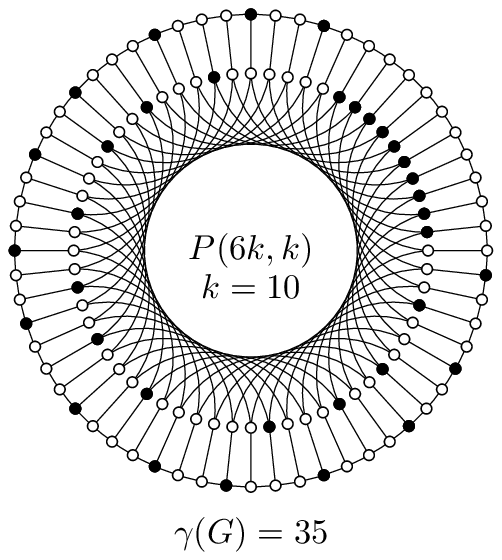} \hspace{2pt}
\includegraphics[scale=0.5]{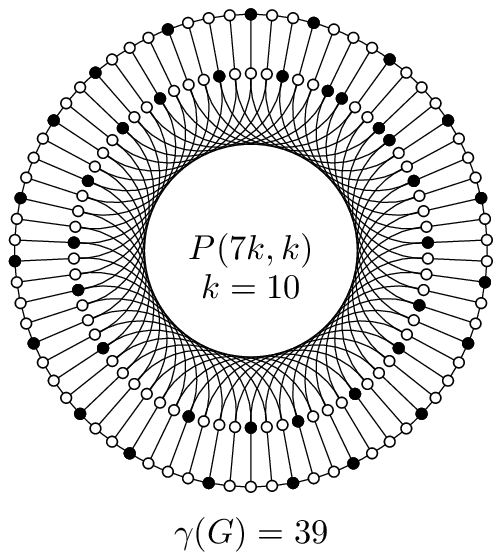} \hspace{2pt}

\centering \vspace{2pt}\small{Figure 2.2:  The dominating sets of
$P(ck,k)$ for $k=4,6,8,10$ and $c=4,5,6,7$}
\end{figure}

\section{The domination number of $P(5k,k)$}

\ \ \ \ In this section, we shall determine the exact domination
number of $P(5k,k)$ for $k\geq 1$.

From Theorem 2.1, we have the following upper bound for $P(5k,k)$.

\noindent \textbf{Lemma 3.1.} For $k\geq 4$, $\gamma(P(5k,k))\leq
3k$.

To prove the lower bound, we need some further notations. In the
rest of the paper, let $S$ be an arbitrary dominating set of
$P(ck,k)$. For convenience, let
$$\begin{array}{llll}
        \ \ \ A_i&=&\{v_{i+jk}: 0\leq j\leq c-1\}, \\
        \ \ \ B_i&=&\{u_{i+jk}: 0\leq j\leq c-1\}, \\
        \ \ \ D_{i(j)}&=&\{v_{i+jk},u_{i+jk}\}, \ \ 0\leq j\leq c-1,\\
\end{array}$$
for $0\leq i\leq k-1$, where the vertices of $A_i$ are on the outer
cycle and those of $B_i$ are on the inner cycle(s). For $0\leq i\leq
k-1$, let $G_i=\langle A_i\cup B_i\rangle$ be the $i$th subgraph
induced by $A_i\cup B_i$ and $S_i=V(G_i)\cap S$.

\noindent \textbf{Lemma 3.2.} Let $\ell\in\{0,1,\ldots,k-1\}$. If
there exists two vertices $v_{x},v_{y} \in  S_{\ell}$ such that
$|x-y|\in \{2k,3k\}$, then $|S_{\ell}|\geq 4$.

\begin{proof} Suppose to the contrary that $|S_{\ell}|\leq 3$.
Without loss of generality, we may assume $x=\ell$ and $y=\ell+2k$,
i.e., $v_{\ell},v_{\ell+2k}\in S_{\ell}$ (see Figure 3.1). Then
at least one vertex of $\{u_{\ell+k},u_{\ell+3k},u_{\ell+4k}\}$
would not be dominated by $S$, a contradiction.
\end{proof}
\begin{figure}[ht]
\centering
\includegraphics[scale=0.85]{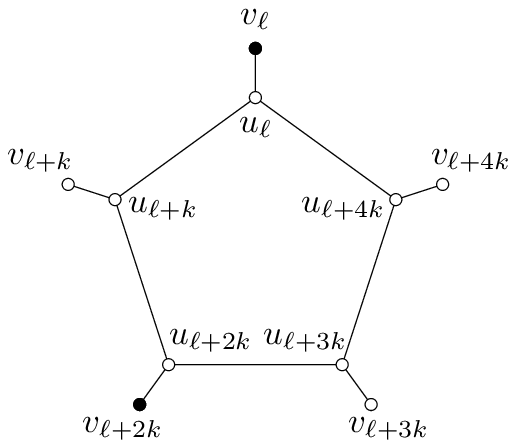}

\centering \vspace{5pt}\small{Figure 3.1:  The graph for the proof
of Lemma 3.1}
\end{figure}

\noindent \textbf{Lemma 3.3.} For any $i\in \{0,1,\ldots,k-1\}$,
$|S_i|\geq 2$. Moreover, if there exists an integer $\ell\in
\{0,1,\ldots,k-1\}$ such that $|S_{\ell}|=2$, then
$S_{\ell}\subseteq B_{\ell}$, $S_{\ell}$ is an independent set, and
the following statements hold.

\indent(i) If $|S_{\ell+1}|=2$, then $|S_{\ell+2}|\geq 4$. Moreover, the equality holds only if $|S_{\ell+3}|\geq 4$; \\
\indent(ii) If $|S_{\ell+1}|=3$, then $|S_{\ell+2}|\geq 3$. Moreover, the equality holds only if $|S_{\ell+3}|\geq 4$; \\
\noindent where the subscripts are taken modulo $k$.

\begin{proof} Since $\langle B_i\rangle$ is isomorphic to $C_5$ and every vertex
of $B_i$ must be dominated by $S_i$, we have that $|S_i|\geq 2$ for
any $i\in \{0,1,\ldots,k-1\}$.

Suppose that there exists an integer $\ell\in \{0,1,\ldots,k-1\}$
such that $|S_{\ell}|=2$.

Assume to the contrary that $|S_{\ell}\cap
B_{\ell}|\leq 1$, or $|S_{\ell}\cap B_{\ell}|=2$ and $S_{\ell}$ is not an independent set.
Then at least one vertex of $B_{\ell}$ would not be dominated by $S$, a contradiction.
Hence, $S_{\ell}\subseteq
B_{\ell}$ and $S_{\ell}$ is an independent set.


(i) Suppose $|S_{\ell+1}|=2$. Then $S_{\ell}\cap
A_{\ell}=\emptyset$, $S_{\ell+1}\cap A_{\ell+1}=\emptyset$ and
$S_{\ell+1}$ is an independent set. Without loss of generality, we
may assume $S_{\ell+1}=\{u_{\ell+1},u_{\ell+1+2k}\}$. Since
$S_{\ell}\cap A_{\ell}=\emptyset$, to dominate
$\{v_{\ell+1+k},v_{\ell+1+3k},v_{\ell+1+4k}\}$, we have
$v_{\ell+2+k},v_{\ell+2+3k},v_{\ell+2+4k}\in S_{\ell+2}$. It follows
from Lemma 3.2 that $S_{\ell+2}\geq 4$.

If $S_{\ell+2}=4$, to dominate $\{u_{\ell+2},u_{\ell+2+2k}\}$, then
$u_{\ell+2+k}\in S_{\ell+2}$, which implies that
$S_{\ell+2}=\{v_{\ell+2+k},v_{\ell+2+3k},v_{\ell+2+4k},u_{\ell+2+k}\}$
and $|D_{\ell+2(0)}\cap S_{\ell+2}|$ $=|D_{\ell+2(2)}\cap
S_{\ell+2}|=0$. Since $S_{\ell+1}\cap A_{\ell+1}=\emptyset$, to
dominate $\{v_{\ell+2},v_{\ell+2+2k}\}$, we have
$v_{\ell+3},v_{\ell+3+2k}\in S_{\ell+3}$ (see Figure 3.2 (1)). It
follows from Lemma 3.2 that $S_{\ell+3}\geq 4$.
\begin{figure}[ht]
\centering
\includegraphics[scale=0.7]{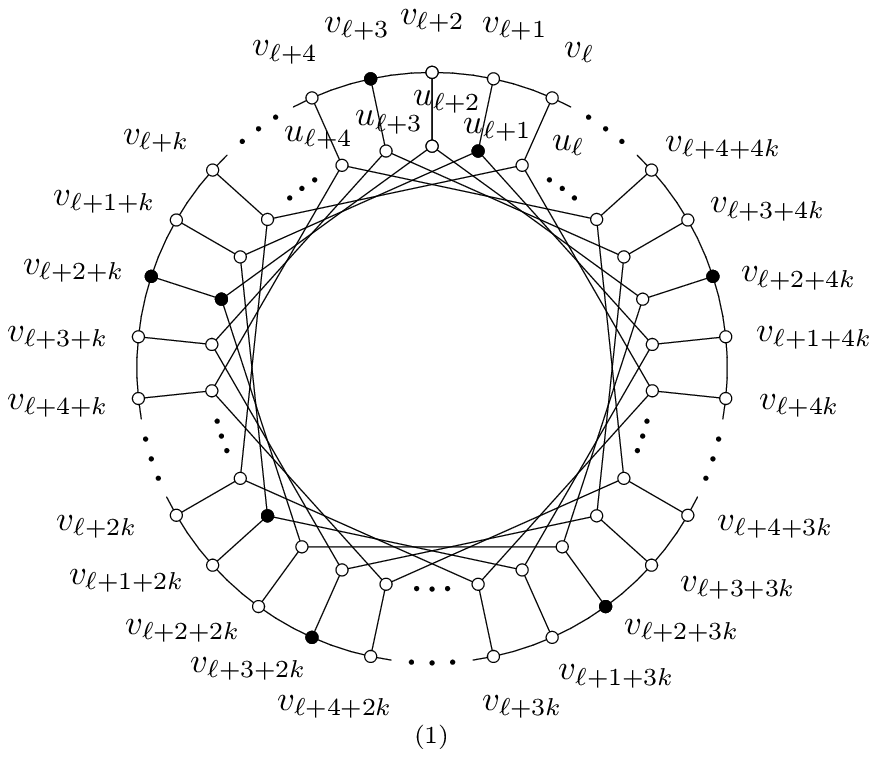} \hspace{10pt}
\includegraphics[scale=0.7]{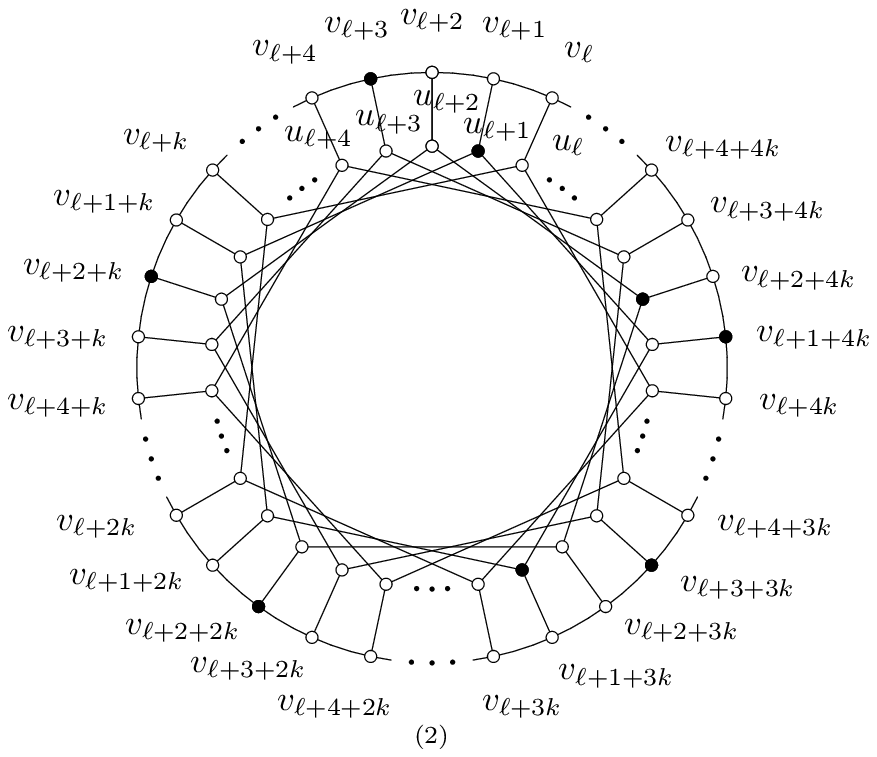}

\centering \vspace{5pt}\small{Figure 3.2:  The graph for the proof
of Lemma 3.2}
\end{figure}

(ii) Suppose $|S_{\ell+1}|=3$. If $|S_{\ell+2}|=2$, then
$S_{\ell}\cap A_{\ell}=\emptyset$ and $S_{\ell+2}\cap
A_{\ell+2}=\emptyset$. To dominate all the vertices in $A_{\ell+1}$,
we have that $|D_{\ell+1(j)}\cap S_{\ell+1}|\geq 1$ for every
$j\in\{0,1,2,3,4\}$. It follows that $|S_{\ell+1}|\geq 5$, a
contradiction with $|S_{\ell+1}|=3$. Hence, $|S_{\ell+2}|\geq 3$.

Now suppose $|S_{\ell+2}|=3$. It is easy to see that there exist at
least two different index $j_1,j_2\in\{0,1,2,3,4\}$ such that
$D_{\ell+1(j_1)}\cap S_{\ell+1}=\emptyset$ and $D_{\ell+1(j_2)}\cap
S_{\ell+1}=\emptyset$.

If $|j_1-j_2|\not\in\{1,4\}$, that is, $|j_1-j_2|\in\{2,3\}$, say
$j_1=1$ and $j_2=3$, since $S_{\ell}\cap A_{\ell}=\emptyset$, to
dominate $\{v_{\ell+1+k},v_{\ell+1+3k}\}$, we have that
$v_{\ell+2+k},v_{\ell+2+3k}\in S_{\ell+2}$. It follows from Lemma
3.2 that $|S_{\ell+2}|\geq 4$, a contradiction with
$|S_{\ell+2}|=3$. Hence, we conclude that $|j_1-j_2|\in\{1,4\}$ and
$|D_{\ell+1(t)}\cap S_{\ell+1}|=1$ for $t\in\{0,1,2,3,4\}\setminus
\{j_1,j_2\}$.

Without loss of generality, we may assume $j_1=1$ and $j_2=2$. To
dominate $\{u_{\ell+1+k}$, $u_{\ell+1+2k}\}$, we have that
$u_{\ell+1},u_{\ell+1+3k}\in S_{\ell+1}$ and $v_{\ell+1}$,
$v_{\ell+1+3k}\not\in S_{\ell+1}$. Since $S_{\ell}\cap
A_{\ell}=\emptyset$, to dominate $\{v_{\ell+1+k},v_{\ell+1+2k}\}$,
we have $v_{\ell+2+k},v_{\ell+2+2k}\in S_{\ell+2}$. Since
$S_{\ell+2}=3$, to dominate $\{u_{\ell+2},u_{\ell+2+3k}\}$, we have
that $u_{\ell+2+4k}\in S_{\ell+2}$. It follows that
$D_{\ell+2(0)}\cap S_{\ell+2}=\emptyset$ and $D_{\ell+2(3)}\cap
S_{\ell+2}=\emptyset$. Since $v_{\ell+1},v_{\ell+1+3k}\not\in
S_{\ell+1}$, we have $v_{\ell+3},v_{\ell+3+3k}\in S_{\ell+3}$ (see
Figure 3.2 (2) for $v_{\ell+1+4k}\in S_{\ell+1}$). It follows from
Lemma 3.2 that $|S_{\ell+3}|\geq 4$.
\end{proof}

\noindent \textbf{Lemma 3.4.} For $k\geq 4$, $\gamma(P(5k,k))\geq
3k$.

\begin{proof} Let $S$ be a dominating set of
$P(5k,k)$ with the minimum cardinality. If $|S_i|\geq 3$ for every
$i\in \{0,1,\ldots,k-1\}$, then
$\gamma(P(5k,k))=|S|=\sum\limits_{i=0}^{k-1}|S_i|\geq 3k$, and we
are done. Hence, we may assume that there exists at least one index
$\ell\in \{0,1,\ldots,k-1\}$ such that $|S_{\ell}|=2$.

Let $H=\{0\leq i\leq n-1: |S_i|=2, |S_{i-1}|>2\}$ and let $h=|H|$.
Let $t_1,t_2,\ldots,t_h$ be all the integers of $H$, where $0\leq
t_1<t_2<\cdots <t_h\leq n-1$. Let $N_i=\{0\leq x\leq n-1: t_i \leq
x\leq t_{i+1}-1\}$ for $i=1,2,\ldots,h$ (In particular,
$t_{h+1}=t_1$). Clearly, $\{0,1,\ldots,n-1\}=\bigcup\limits_{i=1}^h
N_i$. By Lemma 3.3, we conclude that for any $1\leq i\leq h$, $N_i$
satisfies one of the following conditions:

(a) $|S_{t_i}|=2, |S_{t_i+1}|=2, |S_{t_i+2}|\geq 5$ and $|S_x|\geq
3$ for any $t_i+3\leq x \leq t_{i+1}-1$;

(b) $|S_{t_i}|=2, |S_{t_i+1}|=2, |S_{t_i+2}|=4, |S_{t_i+3}|\geq 4$,
$|S_x|\geq 3$ for any $t_i+4\leq x \leq t_{i+1}-1$;

(c) $|S_{t_i}|=2, |S_{t_i+1}|=3, |S_{t_i+2}|\geq 4$, $|S_x|\geq 3$
for any $t_i+3\leq x \leq t_{i+1}-1$;

(d) $|S_{t_i}|=2, |S_{t_i+1}|=3, |S_{t_i+2}|=3, |S_{t_i+3}|\geq 4$,
$|S_x|\geq 3$ for any $t_i+4\leq x \leq t_{i+1}-1$;

(e) $|S_{t_i}|=2, |S_{t_i+1}|\geq 4$, $|S_x|\geq 3$ for any
$t_i+2\leq x \leq t_{i+1}-1$.

It is easy to check that $\sum\limits_{x\in N_i}|S_x|\geq 3|N_i|$
for every $i\in \{1,2,\ldots,h\}$. It follows that
$\gamma(P(5k,k))=|S|=\sum\limits_{0\leq x\leq
k-1}|S_x|=\frac{1}{5}\sum\limits_{0\leq x\leq
n-1}|S_x|=\frac{1}{5}\sum\limits_{i=1}^h\sum\limits_{x\in N_i}
|S_x|\geq \frac{1}{5}\sum\limits_{i=1}^h
3|N_i|=\frac{3}{5}\sum\limits_{i=1}^h |N_i|=\frac{3n}{5}=3k$.
\end{proof}

As an immediate consequence of Lemma 3.1 and Lemma 3.4, we have the
following

\noindent \textbf{Theorem 3.5.} For $k\geq 4$, $\gamma(P(5k,k))=3k$.

\vspace{15pt}

It was shown in \cite{EJM09} that
$\gamma(P(n,1))=\lceil\frac{n}{2}\rceil$ for $n\not\equiv2 \pmod 4$,
$\gamma(P(n,2))$ $=\lceil\frac{3n}{5}\rceil$, and
$\gamma(P(n,3))=\lceil\frac{n}{2}\rceil+1$ for $n\equiv3 \pmod 4$
and $n\neq 11$. Then, we have that $\gamma(P(5,1))=3$,
$\gamma(P(10,2))=6$ and $P(15,3)=9$, which implies that $P(5k,k)=3k$
for $k\in\{1,2,3\}$. Hence, we have the following corollary.

\noindent \textbf{Corollary 3.6.} For $k\geq 1$,
$\gamma(P(5k,k))=3k$.

\section{The domination number of $P(6k,k)$}

\ \ \ \ In this section, we shall determine the exact domination
number of $P(6k,k)$ for $k\geq 1$.

\noindent \textbf{Lemma 4.1.} For $k\geq 4$, $\gamma(P(6k,k))\leq
\lceil\frac{10k}{3}\rceil$.

\begin{proof} To show that $\gamma(P(6k,k))\leq \lceil\frac{10k}{3}\rceil$ for
$k\geq 4$, it suffices to construct a set $S$ that uses
$\lceil\frac{10k}{3}\rceil$ vertices to dominate $P(6k,k)$.

Let $m=\lfloor\frac{k}{3}\rfloor$ and $t=k \mbox{ mod } 3$. Then
$k=3m+t$. Denote {\footnotesize
$$\begin{array}{llll}
S=\left \{\begin{array}{llll}
               \{u_i:0\leq i\leq k-1\}\cup\{u_i:3k\leq i\leq 4k-1\}\cup \\
               \{v_{k+3i+1}:0\leq i\leq\frac{2k}{3}-1\}\cup\{v_{4k+3i+1}:0\leq i\leq\frac{2k}{3}-1\}, & \mbox{ if } t=0;\\
               \\
               \{u_i:0\leq i\leq k-1\}\cup\{u_i:3k-2\leq i\leq 4k-3\}\cup \\
               \{v_{k+3i+1}:0\leq i\leq\frac{2k-2}{3}-1\}\cup\{v_{4k+3i-1}:0\leq i\leq\frac{k-1}{3}-1\}\cup \\
               \{v_{5k-2},v_{5k-1}\}\cup\{v_{5k+3i+2}:0\leq i\leq\frac{k-1}{3}-1\}, & \mbox{ if } t=1;\\
               \\
               \{u_i:0\leq i\leq k-1\}\cup\{v_{k+3i+1}:0\leq i\leq\frac{k-2}{3}\}\cup \\
               \{v_{2k+3i+2}:0\leq i\leq\frac{k-5}{3}-1\}\cup\{u_i:3k\leq i\leq 4k-5\}\cup \\
               \{v_{4k+3i+1}:0\leq i\leq\frac{k-5}{3}-1\}\cup\{v_{5k+3i}:0\leq i\leq\frac{k-2}{3}\}\cup \\
               \{u_{3k-4},v_{3k-2},u_{4k-3},u_{4k-1},v_{4k-3},u_{5k-2},v_{5k-4}\}, & \mbox{ if } t=2.\\
              \end{array}
           \right . \\
\end{array}$$}
It is easy to check that {\footnotesize
$$\begin{array}{llll}
|S|=\left \{\begin{array}{llll}
               2\times3m+2\times\frac{2\times3m}{3}=\lceil\frac{10k}{3}\rceil, & \mbox{ if } t=0;\\
               2\times(3m+1)+\frac{2\times(3m+1)-2}{3}+2\times\frac{3m}{3}+2=\lceil\frac{10k}{3}\rceil, & \mbox{ if } t=1;\\
               2\times(3m+2)-4+2\times(\frac{3m}{3}+1)+2\times\frac{3m-3}{3}+7=\lceil\frac{10k}{3}\rceil, & \mbox{ if } t=2.\\
              \end{array}
           \right . \\
 \end{array}$$}

For $k\equiv 0,1 \pmod 3$, it is not hard to verify that each vertex
in $V(P(6k,k))\setminus S$ can be dominated by $S$.

For $k\equiv 2 \pmod 3$, we have that {\footnotesize
$$\begin{array}{llll}
v_j\in \left \{\begin{array}{llll}
               N[\{u_i:0\leq i\leq k-1\}], & \mbox{ if } 0\leq j\leq k-1;\\
               N[\{v_{k+3i+1}:0\leq i\leq\frac{k-2}{3}\}], & \mbox{ if } k\leq j\leq 2k-1;\\
               N[\{v_{2k+3i+2}:0\leq i\leq\frac{k-5}{3}-1\}\cup\{u_{3k-4},v_{3k-2}\}], & \mbox{ if } 2k\leq j\leq 3k-1;\\
               N[\{u_i:3k\leq i\leq 4k-5\}\cup\{u_{4k-3},u_{4k-1},v_{4k-3}\}], & \mbox{ if } 3k\leq j\leq 4k-1;\\
               N[\{v_{4k+3i+1}:0\leq i\leq\frac{k-5}{3}-1\}\cup\{u_{5k-2},v_{5k-4}\}], & \mbox{ if } 4k\leq j\leq 5k-1;\\
               N[\{v_{5k+3i}:0\leq i\leq\frac{k-2}{3}\}], & \mbox{ if } 5k\leq j\leq 6k-1;\\
              \end{array}
           \right . \\
\end{array}$$}
and {\footnotesize
$$\begin{array}{llll}
u_j\in \left \{\begin{array}{llll}
               N[\{u_i:0\leq i\leq k-1\}], & \mbox{ if } j\in\{\ell k,\ell k+1,\ldots,\ell k+k-1\} \\
               & \mbox{ and } \ell\in\{0,1,5\};\\
               N[\{u_i:3k\leq i\leq 4k-5\}], & \mbox{ if } j\in\{\ell k,\ell k+1,\ldots,\ell k+k-5\} \\
               & \mbox{ and } \ell\in\{2,3,4\};\\
               N[\{u_{3k-4},u_{4k-3},v_{3k-2},u_{4k-1}\}], & \mbox{ if } 3k-4\leq j\leq 3k-1;\\
               N[\{u_{3k-4},u_{4k-3},u_{5k-2},u_{4k-1}\}], & \mbox{ if } 4k-4\leq j\leq 4k-1;\\
               N[\{v_{5k-4},u_{4k-3},u_{5k-2},u_{4k-1}\}], & \mbox{ if } 5k-4\leq j\leq 5k-1.\\
              \end{array}
           \right . \\
\end{array}$$}
\noindent Hence, $S$ is a dominating set of $P(6k,k)$ for $k\geq 4$
with $|S|=\lceil\frac{10k}{3}\rceil$.
\end{proof}

In Figure 4.1, we show the dominating sets of $P(6k,k)$ for $4\leq
k\leq 12$, where the vertices of dominating sets are in dark.
\begin{figure}[ht]
\centering
\includegraphics[scale=0.7]{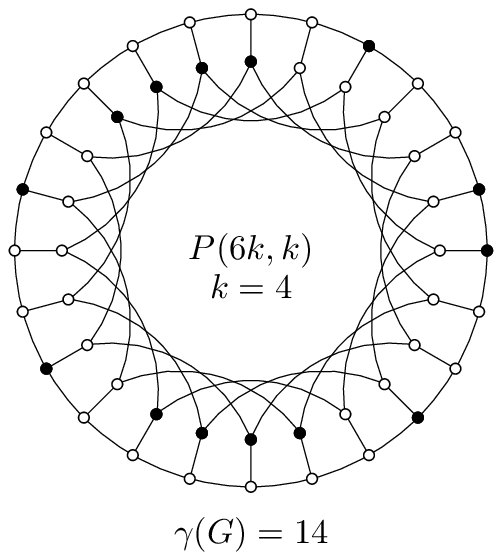}
\includegraphics[scale=0.7]{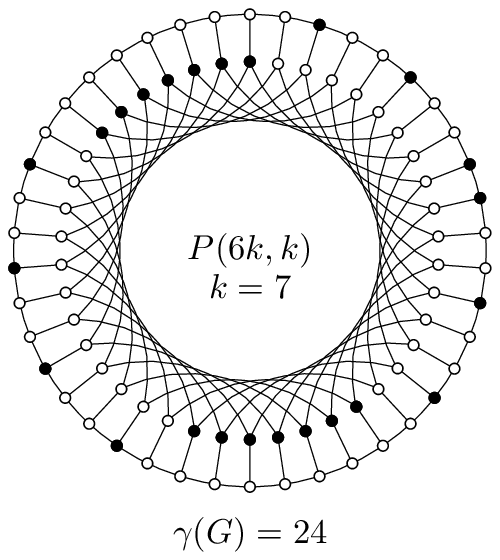}
\includegraphics[scale=0.7]{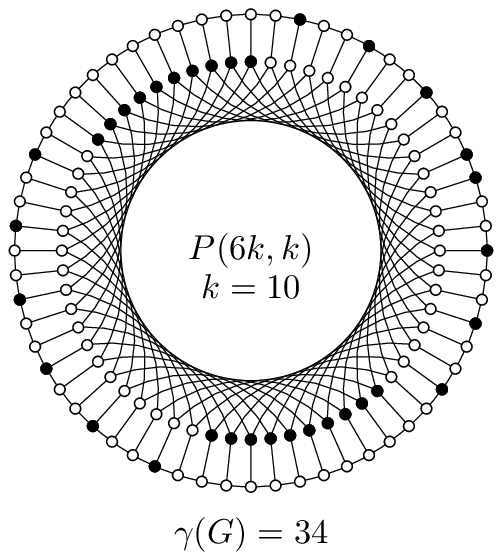}
\includegraphics[scale=0.7]{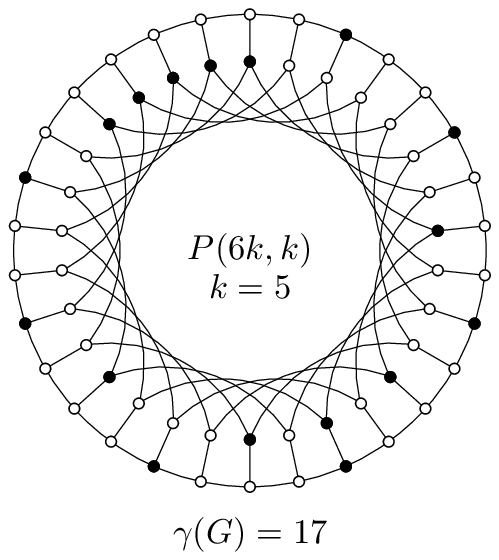}
\includegraphics[scale=0.7]{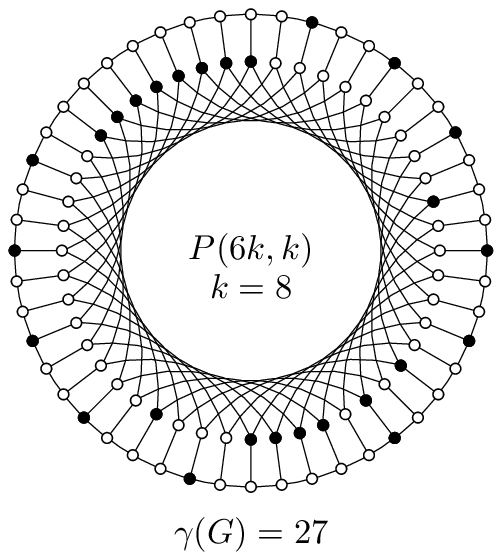}
\includegraphics[scale=0.7]{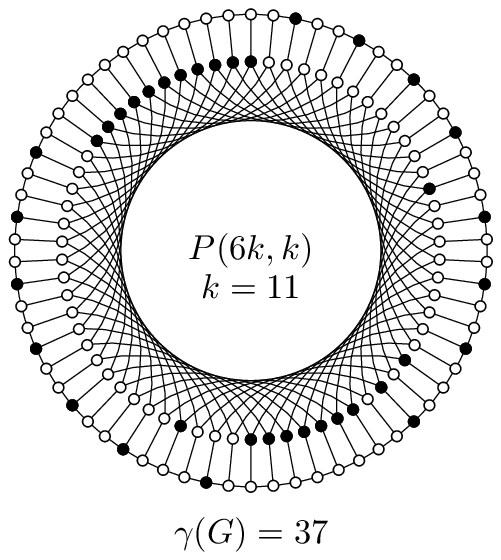}
\includegraphics[scale=0.7]{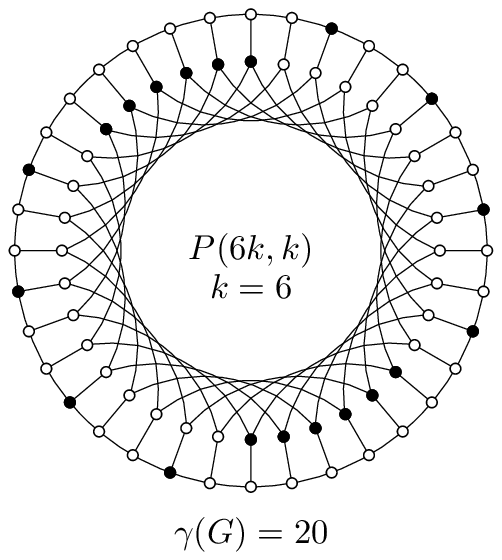}
\includegraphics[scale=0.7]{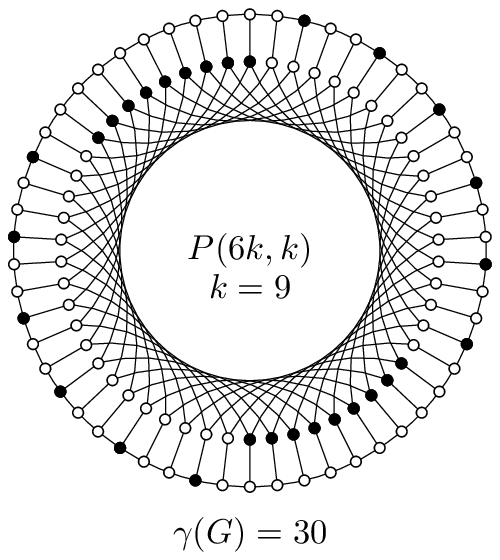}
\includegraphics[scale=0.7]{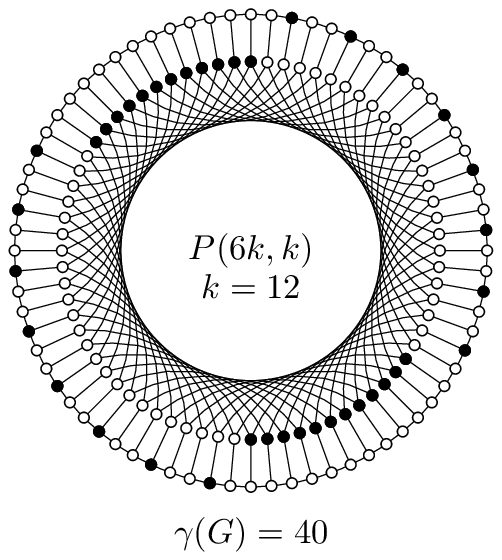}

\centering \vspace{5pt}\small{Figure 4.1:  The dominating sets of
$P(6k,k)$ for $4\leq k\leq 12$}
\end{figure}

\noindent \textbf{Lemma 4.2.} For $i\in \{0,1,\ldots,k-1\}$,
$|S_i|\geq 2$. If there exists an integer $\ell\in \{0,1,\ldots,k-1\}$ such that $|B_{\ell}\cap S_{\ell}|=1$, then $|S_{\ell}|\geq 4$.

\begin{proof} Since $\langle B_i\rangle$ is isomorphic to $C_6$ and every vertex
of $B_i$ must be dominated by $S_i$, we have that $|S_i|\geq 2$ for every $i\in \{0,1,\ldots,k-1\}$.
If there exists an integer $\ell\in \{0,1,\ldots,k-1\}$ such that $|B_{\ell}\cap S_{\ell}|=1$, say $u_{\ell}\in S_{\ell}$,
to dominate $\{u_{\ell+2k},u_{\ell+3k},u_{\ell+4k}\}$, we have $v_{\ell+2k},v_{\ell+3k},v_{\ell+4k}\in S_{\ell}$. It follows that
$|S_{\ell}|\geq 4$. The lemma follows.
\end{proof}

\noindent \textbf{Lemma 4.3.} For every $i\in\{0,1,\ldots,k-1\}$,
$|S_{i-1}\cup S_i\cup S_{i+1}|\geq 10$, where the subscripts are
taken modulo $k$.

\begin{proof} Suppose to the contrary that there exists an integer
$\ell\in\{0,1,\ldots,k-1\}$ such that $|S_{\ell-1}\cup S_{\ell}\cup
S_{\ell+1}|\leq 9$. Combining with Lemma 4.2, we have that
\begin{equation}\label{S_t leq 5}
2\leq |S_t|\leq 5
\end{equation}
for every $t\in\{\ell-1,\ell,\ell+1\}$.

It is easy to see that $V(G_{\ell-1})\cup V(G_{\ell})\cup
V(G_{\ell+1})=(\bigcup\limits_{j=0}^{5}N[v_{\ell+jk}])\cup
B_{\ell-1}\cup B_{\ell+1}$. To dominate each vertex in $A_{\ell}$,
we have that
\begin{equation}\label{N_v geq 1}
|N[v_{\ell+jk}]\cap (S_{\ell-1}\cup S_{\ell}\cup S_{\ell+1})|\geq 1
\end{equation}
\noindent for $0\leq j\leq 5$. It follows that $\sum\limits_{j=0}^{5}|N[v_{\ell+jk}]\cap (S_{\ell-1}\cup S_{\ell}\cup S_{\ell+1})|\geq 6$.
From the assumption, we have
$|(B_{\ell-1}\cap S_{\ell-1})\cup (B_{\ell+1}\cap S_{\ell+1})|\leq
3$. It follows that
\begin{equation}\label{B leq 1}
|B_{\ell-1}\cap S_{\ell-1}|\leq 1 \mbox { \ \ or \ \ } |B_{\ell+1}\cap S_{\ell+1}|\leq 1.
\end{equation}
\indent If $B_{\ell-1}\cap S_{\ell-1}=\emptyset$ or $B_{\ell+1}\cap
S_{\ell+1}=\emptyset$, say $B_{\ell-1}\cap S_{\ell-1}=\emptyset$, to
dominate each vertex in $B_{\ell-1}$, we have $A_{\ell-1}\subset
S_{\ell-1}$, i.e., $|S_{\ell-1}|=6$, a contradiction with \eqref{S_t
leq 5}. Hence,
\begin{equation}\label{B geq 1}
|B_{\ell-1}\cap S_{\ell-1}|\geq 1 \mbox { \ \ and \ \ } |B_{\ell+1}\cap S_{\ell+1}|\geq 1.
\end{equation}
\indent It follows from \eqref{B leq 1} and \eqref{B geq 1} that
$|B_{\ell-1}\cap S_{\ell-1}|=1$ or $|B_{\ell+1}\cap S_{\ell+1}|=1$,
say $|B_{\ell-1}\cap S_{\ell-1}|=1$. Without loss of generality, we
may assume $u_{\ell-1}\in S_{\ell-1}$. To dominate
$\{u_{\ell-1+2k},u_{\ell-1+3k},u_{\ell-1+4k}\}$, we have
$v_{\ell-1+2k},v_{\ell-1+3k}$, $v_{\ell-1+4k}\in S_{\ell-1}$, which
implies $$|S_{\ell-1}|\geq 4.$$ \indent To dominate $u_{\ell+3k}$,
we have that $|\{u_{\ell+2k},u_{\ell+3k},u_{\ell+4k}$,
$v_{\ell+3k}\}\cap S_{\ell}|\geq 1$. It follows that
$\sum\limits_{j=2}^{4}|N[v_{\ell+jk}]\cap (S_{\ell-1}\cup
S_{\ell}\cup S_{\ell+1})|\geq 3+1=4$. Combining with \eqref{N_v geq
1}, we conclude that $|(V(G_{\ell-1})\cup V(G_{\ell})\cup
V(G_{\ell+1})\setminus B_{\ell+1})\cap (S_{\ell-1}\cup S_{\ell}\cup
S_{\ell+1})|=|B_{\ell-1}\cap
S_{\ell-1}|+\sum\limits_{j=0}^{5}|N[v_{\ell+jk}]\cap (S_{\ell-1}\cup
S_{\ell}\cup S_{\ell+1})|\geq 1+7=8$. Hence, we have
$$|B_{\ell+1}\cap S_{\ell+1}|\leq 1.$$ By \eqref{B geq 1}, we have
$|B_{\ell+1}\cap S_{\ell+1}|=1$. It follows from Lemma 4.2 that
$|S_{\ell}|\geq 2$ and $|S_{\ell+1}|\geq 4$. Since $|S_{\ell-1}|\geq
4$, we have $|S_{\ell-1}\cup S_{\ell}\cup S_{\ell+1}|\geq 4+2+4=10$,
a contradiction with assumption. The lemma follows.
\end{proof}

\noindent \textbf{Lemma 4.4.} For $k\geq 4$, $\gamma(P(6k,k))\geq
\lceil\frac{10k}{3}\rceil$.

\begin{proof} Let $S$ be a dominating set of
$P(6k,k)$ with the minimum cardinality. Notice that each subset
$S_i$ is counted 18 times in
$\sum\limits_{i=0}^{6k-1}(|S_i|+|S_{i+1}|+|S_{i+2}|)$. By Lemma 4.3,
we have
$$18\times |S|=\sum\limits_{i=0}^{6k-1}(|S_i|+|S_{i+1}|+|S_{i+2}|)\geq 6k\times
10=60k,$$ which implies that $\gamma(P(6k,k))=|S|\geq
\lceil\frac{10k}{3}\rceil$.
\end{proof}

As an immediate consequence of Lemma 4.1 and Lemma 4.4, we have the
following

\noindent \textbf{Theorem 4.5.} For $k\geq 4$, $\gamma(P(5k,k))=3k$.

\vspace{15pt}

It was shown in \cite{EJM09} that $\gamma(P(n,1))=\frac{n}{2}+1$ for
$n\equiv2 \pmod 4$, $\gamma(P(n,2))$ $=\lceil\frac{3n}{5}\rceil$ and
$\gamma(P(n,3))=\frac{n}{2}+1$ for $n\equiv2 \pmod 4$. Then, we have
that $\gamma(P(6,1))=4$, $\gamma(P(12,2))=8$ and $P(18,3)=10$, which
implies that $P(6k,k)=\lceil\frac{10k}{3}\rceil$ for $k\in\{1,3\}$
and $P(6k,k)=\lceil\frac{10k}{3}\rceil+1$ for $k=2$. Hence, we have
the following corollary.

\noindent \textbf{Corollary 4.6.} For $k\geq 1$,
$$\begin{array}{llll}
\gamma(P(6k,k))=& \left \{\begin{array}{llll}
               \lceil\frac{10k}{3}\rceil, & \mbox{ if } k\neq2;\\
               \lceil\frac{10k}{3}\rceil+1, & \mbox{ if } k=2.\\
              \end{array}
           \right . \\
\end{array}$$

\end{document}